\documentclass[12pt]{amsart}
\usepackage{geometry}                
\usepackage{graphicx}
\usepackage{amsthm}
\usepackage{amssymb}
\usepackage{epstopdf}
\usepackage{float}

\usepackage{fancyhdr}
\newcommand{\lp}{\left(}
\newcommand{\rp}{\right)}
\newcommand{\lb}{\left[}
\newcommand{\rb}{\right]}
\newcommand{\lc}{\left\{}
\newcommand{\rc}{\right\}}
\newcommand{\ddx}[1]{\operatorname{\frac{d#1}{dx}}}
\newcommand{\ddxk}[1]{{\frac{\operatorname{d}^{#1}}{\operatorname{dx}^{#1}}}}
\newcommand{\x}{\operatorname{x}}
\newcommand{\summ}[3]{\sum_{#1=#2}^{#3}}

\newcommand{\R}{\mathbb{R}}
\newcommand{\rank}{\operatorname{rank}}

\theoremstyle{definition}
\newtheorem{definition}{Definition}
\newtheorem{exmp}{Example}[section]

\theoremstyle{plain}

\newtheorem{theorem}{Theorem}
\newtheorem{corollary}{Corollary}

\title{Discrete approximations of differential equations via trigonometric interpolation}
\author{Oksana Bihun${}^1$}
\thanks{${}^1$Corresponding author's e-mail: obihun@cord.edu.}
\author{Austin Bren${}^2$} 
\thanks{${}^2$E-mail: asbren@cord.edu} 
\author{Michael Dyrud${}^3$}
\thanks{${}^3$E-mail: mdyrud@cord.edu} 
\author{Kristin Heysse${}^4$}
\thanks{${}^4$E-mail: keheysse@cord.edu}
\date{\today}  
\begin{document}

\maketitle

\vspace{-7mm}
\begin{center}
\small
{Concordia College, 901 8th Street South, Moorhead, MN 56560, USA}
\normalsize
\end{center}

\begin{abstract}
To approximate solutions of a linear differential equation, we project,  via trigonometric interpolation, its solution space onto a finite-dimensional space of trigonometric polynomials and construct a matrix representation of the differential operator associated with the equation.  We compute the ranks of the matrix representations of a certain class of linear differential operators.  Our numerical tests show high accuracy and fast convergence of the method applied to several boundary and eigenvalue problems.
\end{abstract}

\textbf{MSC 34K28, 65N22, 15A03, 34B05, 34B09, PACS 0002}
\section{Introduction}
In this paper,  we use trigonometric interpolation to approximate solutions of a differential equation $\mathcal{A}u=f$, whose differential operator $\mathcal{A}$ with domain $D(\mathcal{A})$ is a formal polynomial of operators $\lc1,\x,\ddx{} \rc$, and $f \in \operatorname{Range}(\mathcal{A})$.   A solution $u$ is projected onto the space $\mathcal{T}_n^L$ of $L$-periodic trigonometric polynomials of degree~$n$. The projection $T[u]$, defined as a trigonometric interpolant of $u$, is identified with a vector $\hat{\mathbf{u}}$ of its values at $N=2n+1$ partition points via an isomorphism $\psi$, where $n \in \mathbb{N}$. The operator $\mathcal{A}$ is represented by an $N \times N$ matrix $A$ defined implicitly by $A\hat{\mathbf{v}}=\psi T\lb \mathcal{A}T[v]\rb$ for all $v \in D(\mathcal{A})$. The original equation is approximated by a system of linear equations $A\hat{\mathbf{u}}=\hat{\mathbf{f}}$, where $\hat{\mathbf{f}}=\psi T [f]$.

In sect. \ref{RankSection}, we prove that if $\mathcal{A}=\alpha_0+\alpha_1\ddx{}+\ldots+\alpha_s\ddxk{s}$, then $$\mbox{rank}\, A=N+|\mbox{sign } \alpha_0|-2m-1,$$ where $m \geq 0$ is the number of solutions, in the set $\{1,2,\ldots,n\}$, of a system of polynomial equations dependent on the set of coefficients $\lc \alpha_i \rc_{i=0}^s$. 

Our method employs ideas similar to those of the method of discrete approximations of differential operators that uses Lagrange interpolation.  The latter approach was developed and tested by Calogero~\cite{Calogero,Calogero2} to solve eigenvalue problems for ordinary and partial differential operators.   Using this method, Mytropolski, Prykarpatsky, and Samoylenko created a general algebraic-projection framework and furthered the applications to Mathematical Physics~\cite{Myt}. 

Using Lagrange interpolation by algebraic polynomials, Calogero introduced a finite-dimensional matrix representation $K \in R^{N \times N}$ of a differential operator $\mathcal{K}$ that is a formal polynomial of $\lc 1,\x,\ddx{} \rc $.  A differential equation $\mathcal{K}u=f$, where $f \in \operatorname{Range}(\mathcal{K})$, is then approximated by a system of linear equations $K\hat{\mathbf{u}}=\hat{\mathbf{f}}$, where $\hat{\mathbf{f}} \in \R^{N}$ is a finite representation of $f$. Bihun and Lu\'{s}tyk estimated the approximation error of the method~\cite{Bihun2,Bihun3}.
 
Bihun and Prytula proved that the rank of the representation $K \in \R^{N \times N}$ of a differential operator $$\mathcal{K}= a_k \ddxk{k}+a_{k+1} \ddxk{k+1}+\ldots+a_{k+p} \ddxk{k+p},$$ is $N-k$ \cite{Bihun}.  Therefore, completion of the discrete approximation $K\hat{\mathbf{u}}=\hat{\mathbf{f}}$ with initial or boundary conditions may lead to an inconsistent system of linear equations.  In~\cite{Bihun,Calogero}, this obstacle was overcome by employing substitutions.  Also, the authors recognized that Lagrange approximation is less effective for approximating periodic solutions of differential equations.  Our study of a discrete approximation $A \in \R^{N\times N}$ of $\mathcal{K}$, obtained using trigonometric interpolation, aims to provide a class of differential operators whose representations' ranks are compatible with the dimension of the solution space of the differential equation $\mathcal{K}u=f$.   

In two short communications~\cite{Calogero3,Calogero4}, Calogero explored the idea of approximation of differential operators using trigonometric interpolation.  He proposed two matrix representations of the differential operator $\ddx{}$ and calculated their eigenvalues.  It was proven that if a trigonometric polynomial $u$ of degree $n$ solves a differential equation $\mathcal{A}u=0$, where $\mathcal{A}=\summ{j}{0}{s} \alpha_j \ddxk{j}$, then the solution of the approximate problem $A\hat{\mathbf{u}}=\mathbf{0}$ can be used to recover the \textit{exact} solution $u$, where $A \in \R^{N\times N}$ is a matrix representation of $\mathcal{A}$ constructed using trigonometric interpolation (\textit{cf.} theorem~\ref{calogero} in subsection~\ref{eigensubsection}).  To our knowledge, the method of discrete approximation of differential equations via trigonometric interpolation has not been implemented and tested.  Our numerical tests show high accuracy and fast convergence of the method applied to several boundary and eigenvalue problems for differential equations.

\section{Trigonometric interpolation}
\begin{definition}
A \textit{trigonometric polynomial} $T(x)$ of period $L$ and degree $n$ is defined to be
\begin{equation}
\label{Polynomial}
T(x)=c+\sum_{k=1}^n {\left(a_k \sin\lp \frac{2\pi}{L}k x\rp + b_k \cos\lp \frac{2\pi}{L}k x\rp\right)},
\end{equation}
where $ c, a_1,...,a_n,b_1,...,b_n $ are real numbers and $a_n^2+b_n^2 > 0$.
\end{definition}

We will refer to trigonometric polynomials as simply polynomials. Should the need arise to discuss ordinary polynomials, they will be referred to as \textit{algebraic polynomials}.  We denote the space of all $L$-periodic trigonometric polynomials of degree $n$ by $\mathcal{T}_{n}^L $.  

In \cite{Bihun,Calogero,Myt}, Lagrange interpolation is used to construct finite-dimensional representations of operators $\left\{ 1, \x, \ddx{} \right\} $.  We construct similar representations using trigonometric polynomials to solve differential equations and eigenvalue problems.  This modification of Calogero's method is especially useful for solving differential equations and eigenvalue problems with periodic solutions. Our motivation to use trigonometric polynomials to approximate solutions of differential equations is also based on the Second Weierstrass Theorem, which states that every continuous $L$-periodic function $f : [0,L] \rightarrow \R$ can be approximated, with any desired accuracy in the supremum norm, by an $L$-periodic trigonometric polynomial \cite{Akhiezer}.

To construct representations of differential operators, we use trigonometric interpolation.  Let $ u : [0,L) \rightarrow \mathbb{R} $ be a continuous function.  Partition the interval $[0,L)$ into $2n+1$ points $\{x_0,x_1,...x_{2n}\}$ such that $0=x_0 < x_1< ... < x_{2n} <L$. Let $u_j = u(x_j)$ for all $j\in \{0,1,...,2n\}$.  We construct a polynomial $T[u] \in \mathcal{T}_n^L$ such that 
\begin{equation}
\label{criteria}
T[u](x_j)=u_j \mbox{ for all } 0 \leq j \leq 2n.
\end{equation}

This polynomial is defined to be~\cite{Zygmund}
\begin{equation}
\label{T}
T[u](x)=\sum_{j=0}^{2n} {u_j t_j (x)},
\end{equation}
where
\begin{equation}
\label{t}
t_j(x)=\frac{\displaystyle\prod_{\substack{k=0\\ k\neq j}}^{2n} { \sin\left(\frac{\pi(x-x_k)}{L}\right)}}{\displaystyle\prod_{\substack{k=0\\ k\neq j}}^{2n} { \sin\left(\frac{\pi(x_j-x_k)}{L}\right)}}
\end{equation}
for all $x \in [0,L).$

The functions $t_j$ satisfy $t_j(x_k) = \delta_{jk}$, where $\delta$ is the Kronecker delta. Moreover, each $t_j$ is a trigonometric polynomial of degree $n$, which follows from trigonometric product-to-sum formulas.  

We note that the set of functions $\lc t_j \rc_{j=0}^{2n}$ is linearly independent.  Indeed, suppose that $\summ{i}{0}{2n} \alpha_i t_i (x)=0$ for all $x \in [0,L)$ and some real numbers $\alpha_0,\alpha_1,\ldots,\alpha_{2n}$. By substituting $x=x_j$ into the last equation, we obtain $\alpha_j=0$ for all $j\in\{0,1,\ldots,2n\}$.

\begin{theorem}
\label{trigUnique}
For a given set of points $S=\lc \lp x_0,y_0 \rp, \lp x_1,y_1 \rp, ...,  \lp x_{2n},y_{2n} \rp \rc \in \mathbb{R}^2$, such that $ x_0,...,x_{2n} \in \lb 0,L \rp$ are distinct, there exists a unique $L$-periodic trigonometric polynomial $f$ of degree $n$ that passes through all the points in $S$.  
\end{theorem}

The existence of a trigonometric polynomial $$f(x)=c+\sum_{k=1}^n\lp a_k\cos\lp\frac{2\pi}{L}kx\rp+ b_k\sin\lp\frac{2\pi}{L}kx\rp \rp$$ whose graph contains the set $S$ follows from the previous discussion, and its uniqueness is a consequence of the Fundamental Theorem of Algebra applied to an algebraic polynomial $h_f: \mathbb{C} \rightarrow \mathbb{C}$ associated with $f$.  More precisely, $f$ is recast as $f(x)=\sum_{m=-n} ^n c_m e^{\frac{2\pi}{L}imx}$, where $i^2=-1$, $c_m=\frac{1}{2}(a_m-ib_m)$, $c_{-m}=\frac{1}{2}(a_m+ib_m)$ for $m \in \{1,2,\ldots,n\}$, and $c_0=c$. The polynomial $h_f$ is defined by $h(z)=z^n \sum_{m=-n}^n c_m z^m$, where the complex variable $z=e^{\frac{2\pi}{L}ix}$. If $g$ is another trigonometric polynomial that passes through the points in $S$, then the algebraic polynomial $Q=h_f-h_g$ of degree at most $2n$ has $2n+1$ distinct zeros $z_k=e^{\frac{2\pi}{L}ix_k}$, $k=0,1,\ldots,2n$. By the Fundamental Theorem of Algebra, $Q=0$, which implies $f=g$.

\begin{corollary}
\label{TofFisF}
Let $f$ be an $L$-periodic trigonometric polynomial of degree $n$, and let the projection $T[f]$ be constructed using a partition $\lc x_i\rc_{i=0}^{2n}$ of $[0,L)$, according to formula~\eqref{T}. Then $T[f]=f$.  Moreover, the system of functions $\{t_j\}_{j=0}^{2n}$ forms a basis of $\mathcal{T}^L_n.$
\end{corollary}

It is known that the set of functions $$\left\{1,\cos\lp\frac{2\pi}{L}x\rp,\sin\lp\frac{2\pi}{L}x\rp,\ldots,\cos\lp n\frac{2\pi}{L}x\rp,\sin\lp n \frac{2\pi}{L}x\rp\right\}$$ forms a basis of $\mathcal{T}_n^L$.  We denote this system by $\{f_i \}_{i=0}^{2n},$ where $f_0(x)=1,f_1(x)=\cos\lp\frac{2\pi}{L}x\rp,f_2(x)=\sin\lp\frac{2\pi}{L}x\rp,\ldots,f_{2n}(x)=\sin\lp\frac{2\pi}{L}nx\rp.$  We will use both bases in our discussion.  

\section{Matrix representations of differential operators via trigonometric interpolation}

In this section we construct finite-dimensional representations of linear differential operators. We denote the space of all $s$ times continuously differentiable functions on $[0,L]$ by $C^s_{[0,L]}$.

\begin{theorem}
\label{A}
Let $\mathcal{A}$ be a linear differential operator of order $s$ with a domain $D(\mathcal{A})\supset C^s_{[0,L]}$. Let $u \in D(\mathcal{A}) $ and let $ \{x_0,x_1,...,x_{2n} \}$ be a partition of $[0,L)$.  Let $\hat{\mathbf{u}}= (u_0,u_1,...,u_{2n})^T$, where $\hat{\mathbf{u}}_j=u(x_j)\mbox{ for all } j \in \{0,1,...,2n\}$.  There exists a $(2n+1)$-dimensional square matrix $A$ such that
\begin{equation}
\label{AEqn}
T\left[\mathcal{A} T[u]\right](x) = \sum_{j=0}^{2n} { [A\hat{\mathbf{u}}]_j t_j(x) }, \qquad x \in [0,L).
\end{equation}
Moreover, 
\begin{equation}
\label{newA}
A_{jk}=\mathcal{A}t_k(x_j)
\end{equation}
for all $j,k \in \{0,1,\ldots,2n\}$.
\begin{proof}
To find a formula for the components of $A$, we substitute $u=t_k$ into eq.~\eqref{AEqn}, where $k \in \{0,1,\ldots,2n\}$.  From Corollary \ref{TofFisF},  $T[t_k]=t_k$ and eq.~\eqref{AEqn} becomes $$T\left[\mathcal{A}t_k\right](x)=\sum_{j=0}^{2n} { \left[A \mathbf{e}_k\right]_j t_j(x) }, \qquad x \in [0,L),$$ 
where $\left\{\mathbf{e}_k\right\}_{k=0}^{2n}$ is the standard basis of $\mathbb{R}^{2n+1}$.

By definition~\eqref{T},  $$\sum_{j=0}^{2n} { \lb \mathcal{A} t_k(x_j) \rb t_j(x) }= \sum_{j=0}^{2n} {A_{jk} t_j(x) }, \qquad x \in [0,L).$$ Therefore,
$\mathcal{A}t_{k}(x_j) = A_{jk} $
for all $j,k \in \{ 0,1,\ldots,2n\}$.
Let us verify that eq.~\eqref{AEqn} holds for an arbitrary function~$u \in D(\mathcal{A})$.  With $\hat{\mathbf{u}}$ defined as in the hypothesis, we obtain   
\begin{align*}
T\left[\mathcal{A}T[u]\right](x) =& T \lb \mathcal{A}\sum_{k=0}^{2n}u_kt_k \rb(x) \\
=& \sum_{j=0}^{2n}\lb\sum_{k=0}^{2n} u_k\mathcal{A}t_k(x_j)\rb t_j(x) \\
=&\sum_{j=0}^{2n} \left[  \sum_{k=0}^{2n} A_{jk} u_{k} \right] t_j(x) \\
=& \sum_{j=0}^{2n} \left[ A\hat{\mathbf{u}}\right]_j t_j(x), \;\;x \in [0,L), 
\end{align*} 
as required.
\end{proof}
\end{theorem}

\begin{definition} If a matrix $A\in\mathbb{R}^{(2n+1)\times(2n+1)}$ and a linear differential operator $\mathcal{A}$ with domain $D(\mathcal{A}) \supset C_{[0,L]}^s$ are related by eq.~\eqref{AEqn}, we say that the matrix $A$ \textit{represents} the differential operator $\mathcal{A}$ in $\mathbb{R}^{(2n+1)}$. \end{definition}

Let $D$ be a $(2n+1)$-dimensional square matrix that represents the operator $\ddx{}$ and let $X$ be a $(2n+1)$-dimensional square matrix that represents the operator of multiplication by $x$. Using theorem~\ref{A}, let us find explicit formulas for the entries of matrices $D$ and $X$. For convenience, we define
\begin{equation}
\label{piconstant}
\psi_m = \displaystyle\prod_{\substack{k=0\\ k\neq m}}^{2n} { \sin\left(\frac{\pi(x_m-x_k)}{ L}\right)}
\end{equation}
for all $m \in \{0,1,\ldots,2n\}$.  For every vector $v\in\mathbb{R}^N$, $\operatorname{Diag}(v)$ denotes the $N$-dimensional square matrix whose diagonal is $v$; $I_N$ denotes the $N$-dimensional identity matrix.

\begin{theorem}
\label{ExplicitDX}

Let $D,X$ be $(2n+1)$-dimensional square matrices that represent the operator $\ddx{}$ and the operator $\x$ of multiplication by  $x$, respectively. Then
\begin{equation}
\label{D}
D_{mj} = \lc \begin{array}{lr} \displaystyle \frac{\psi_m}{\psi_j} \frac{\pi}{L \sin{\frac{\pi (x_m-x_j)}{L}}} & \mbox{ if } j\neq m, \vspace{0.1 in}  \\ \displaystyle \frac{\pi}{L} \sum_{\substack{k=0 \\ k\neq j}}^{2n}  \cot{\frac{\pi(x_j-x_k)}{L}} & \mbox{ if } j=m  \end{array} \right.
\end{equation}
and $X=\operatorname{Diag}(x_0,x_1,\ldots,x_{2n}).$ Moreover, if $P(x)=\alpha_0+\alpha_1x+\ldots+\alpha_sx^s$ is an algebraic polynomial of degree $s$ with real coefficients $\{\alpha_i\}_{i=0}^s$, then the operator $P\lp\ddx{}\rp$ is represented by the matrix $P\lp D\rp$.
\end{theorem}
The proof of this theorem is a straightforward application of formulas~\eqref{AEqn} and~\eqref{newA}.
We note that if $g$ is a trigonometric polynomial of degree $n$, then $\ddx{g}$ is also a trigonometric polynomial of degree $n$.  The representation of the identity operator 1 is easily verified to be $I_{2n+1}$.  

\section{Rank of matrix representations of differential operators}
\label{RankSection}
Let $P(x)=\alpha_0+\alpha_1x+\ldots+\alpha_s x^s$ be a polynomial of degree $s$ and let $\mathcal{A}=P(\ddx{})$ be a differential operator with domain $D(\mathcal{A}) \supset C_{[0,L]}^s$.  By theorem~\ref{ExplicitDX}, the matrix $P(D)$ represents the differential operator $\mathcal{A}$ on $\R^{2n+1}$.  The main idea of our method is to approximate differential equation $\mathcal{A}u=f$, where $f \in C_{[0,L]}$, by the system of linear equations $P(D)\hat{\mathbf{u}}=\hat{\mathbf{f}},$ where $\hat{\mathbf{f}}=\lp f(x_0),f(x_1),\ldots,f(x_{2n}) \rp^T$.  To apply this idea to solution of boundary and initial value problems, we need to study the rank of the matrix $P(D)$.  The following theorem gives a complete answer to this question.  

%Begin Theorem.
\begin{theorem}
\label{DtotheN}
Let $P$ be a polynomial of degree $s$ given by 
\begin{equation}
P(x)=\alpha_0+\alpha_1x+\ldots+\alpha_sx^s,
\end{equation} where $\alpha_0,\alpha_1,\ldots,\alpha_s \in \R$ and $\alpha_{s}\neq 0$. Define the polynomials 
$$\phi(k)=\alpha_0-\alpha_2 \lp \frac{2\pi}{L}k \rp^2+\alpha_4 \lp \frac{2\pi}{L}k \rp^4-\ldots+(-1)^\ell \alpha_q\lp\frac{2\pi}{L}k\rp^q,$$ 
and
\begin{equation*}
\zeta{(k)}=\alpha_1-\alpha_3\lp\frac{2\pi}{L}k\rp^2+\alpha_5\lp\frac{2\pi}{L}k\rp^4+\ldots+(-1)^m \alpha_{r} \lp \frac{2\pi}{L}k \rp^{(r-1)},
\end{equation*}
where $q=2 \lfloor \frac{s}{2} \rfloor$ is the largest even integer that does not exceed $s$, $r=2s-q-1$ is the largest odd integer that does not exceed $s$, $2\ell=q \bmod 4$, and $2 m=(r-1)\bmod 4$.

Let $D$ be the $(2n+1)$-dimensional square matrix, defined by \eqref{D}, that represents the operator $\ddx{}$ in $\mathbb{R}^{2n+1}$. Suppose that the system
\begin{equation}
\label{phizeta}
\Bigg\{ \begin{array}{lr} \phi(k)=0, \\ \zeta(k)=0 \end{array}
\end{equation}
has exactly $m\geq 0$ solutions that belong to the set $\{ 1,2,\ldots,n\}$.  Then 
\begin{equation}
\operatorname{rank}P(D)=2n + |\operatorname{sign } \alpha_0|-2m
\end{equation}

\end{theorem}
\textit{Note.} The results stated in this theorem do not depend on the choice of partition $\{ x_i\}_{i=0}^{2n}$ of $[0,L)$ with which the matrix $D$ is constructed.

\begin{proof}
We prove the theorem for the case $\alpha_0 \neq 0$.  The case $\alpha_0=0$ can be proven in a similar fashion.  

By the rank-nullity theorem, it suffices to show that the dimension of the kernel of $P(D)$ is $2m$.  
Define the differential operator $\mathcal{A} =P(\ddx{})$ and 
let $\hat{\mathbf{u}} \in \R^{2n+1}$.  We observe that $P(D)\hat{\mathbf{u}}=\mathbf{0}$ if and only if $\summ{i}{0}{2n} \lb P(D) \hat{\mathbf{u}} \rb_i{t}_i(x)=0$, which happens if and only if $\mathcal{A}\summ{i}{0}{2n} \hat{\mathbf{u}}_i t_i(x)=0$, where $ x \in [0,L)$. Therefore, we will study the kernel of the differential operator $\mathcal{A}$ with its domain restricted to $\mathcal{T}_n^L$.

Let $g\in \mathcal{T}_n^L$ be given by $$g(x)=b+ \sum_{k=1}^{n} \lp c_k\cos \lp \frac{2\pi}{L}kx\rp+d_k\sin \lp\frac{2\pi}{L}kx\rp \rp, \qquad x \in [0,L). $$ We will prove that $\mathcal{A}g=0$ if and only if $b=0$ and $c_k=d_k=0$ for all $k \in \{1,2,\ldots,n\}$ that do not solve system~\eqref{phizeta}. That is, $\ker \mathcal{A}|_{\mathcal{T}_n^L}=\{0\}$ if $m=0$ and $\ker \mathcal{A}|_{\mathcal{T}_n^L}=\lc \cos\lp\frac{2\pi}{L}k_j x\rp, \sin\lp\frac{2\pi}{L}k_j x\rp \rc_{j=1}^m$ if $m\geq 1$ and $k_1,k_2,\ldots, k_m$ are all the solutions of system~\eqref{phizeta} in $\{1,2,\ldots,n\}$.

Indeed, because
\begin{align}
\label{Q}
\mathcal{A}g(x)=&b\alpha_0+\sum_{k=1}^{n} \left[ c_k \sum_{i=0}^{s} \alpha_i \ddxk{i}  \cos\lp\frac{2\pi}{L}kx\rp  + d_k \sum_{i=0}^{s} \alpha_i \ddxk{i}  \sin\lp\frac{2\pi}{L}kx\rp  \right] \nonumber \\
=&b\alpha_0 + \sum_{k=1}^{n} \left[ \lp \phi(k)c_k + \frac{2\pi k}{L} d_k \zeta(k) \rp \cos\lp\frac{2\pi}{L}kx\rp \right. \nonumber \\ & \qquad + \left. \lp -\frac{2\pi k}{L}  c_k\zeta(k)+d_k\phi(k) \rp \sin\lp\frac{2\pi}{L}kx\rp \right], \qquad x \in [0,L),
\end{align}
$\mathcal{A}g=0$ implies $b=0$,
\begin{align}
\label{firstZeroEq}
c_k \phi(k)+ \frac{2\pi k}{L} d_k \zeta(k)=0, \mbox{ and}
\end{align}
\begin{align}
\label{secondZeroEq}
-\frac{2\pi k}{L} c_k \zeta(k)+ d_k \phi(k) =0
\end{align}
for all $k \in \{1,2,\ldots,n\}$. Suppose that $k \in \{1,2,\ldots,n\}$ does not solve system~\eqref{phizeta}. Then $\phi(k)$ and $\zeta(k)$ are not both equal to zero, and
\begin{equation}
\label{dettilde}
\operatorname{det} \left( {\begin{array}{cc} \phi(k) &  \frac{2\pi k}{L} \zeta(k) \\ - \frac{2\pi k}{L} \zeta(k) & \phi(k) \\ \end{array} } \right)=\phi^2(k)+\lp\frac{2\pi k}{L}\rp^2\zeta^2(k)>0. 
\end{equation}
  Therefore, the linear system (\ref{firstZeroEq},\ref{secondZeroEq}) has the unique solution $c_k=d_k=0.$
  
  On the other hand, if $g$ is a linear combination of $\lc \cos\lp\frac{2\pi}{L}k_j x\rp, \sin\lp\frac{2\pi}{L}k_j x\rp \rc_{j=1}^m$ or $g=0$, then $\mathcal{A}g=0$ by formula~\eqref{Q}.
  
  By a previous observation, the vector $\hat{\mathbf{u}} \in \ker P(D)$ if and only if $\mathcal{A}\sum_{i=0}^{2n}\hat{\mathbf{u}}_i t_i(x)=0$ for all $x \in [0,L)$. If $m=0$, then system~\eqref{phizeta} has no solutions in $\{1,2,\ldots,n\}$, which implies $\sum_{i=0}^{2n}\hat{\mathbf{u}}_i t_i=0$ and $\hat{\mathbf{u}}=\mathbf{0}$.  In this case $\ker P(D)=\{0\}$ and $\rank P(D)$ is full.

If $m \geq 1$ and the numbers $k_1,k_2,\ldots,k_m $ are all the solutions of system  \eqref{phizeta} in $\{ 1,2,\ldots,n \}$,  then 
 $\summ{i}{0}{2n} \hat{\mathbf{u}}_i t_i(x)$ is a linear combination of $\lc \cos \lp \frac{2\pi}{L}k_j x\rp, \sin\lp \frac{2\pi}{L}k_j x\rp \rc_{j=1}^m$.  This occurs if and only if $\hat{\mathbf{u}}$ is a linear combination of the vectors $\lc \mathbf{v}_j, \mathbf{w}_j \rc_{j=1}^m \subset \R^{2n+1}$ defined by the equations $$\cos\lp \frac{2\pi}{L}k_j x\rp =\summ{i}{0}{2n} \lb \mathbf{v}_j \rb_i t_i(x)$$ and  $$\sin\lp \frac{2\pi}{L}k_j x\rp =\summ{i}{0}{2n} \lb \mathbf{w}_j \rb_i t_i(x)$$ for all $j \in \{1,2,\ldots,m\}$. The system of vectors $\lc \mathbf{v}_j, \mathbf{w}_j \rc_{j=1}^m$ is linearly independent because the system of functions $\left\lbrace \cos \lp \frac{2\pi}{L}k_jx\rp,\sin \lp \frac{2\pi}{L}k_jx \rp \right\rbrace_{j=1}^m$ is linearly independent.  Therefore, $\operatorname{ker}P(D)=\operatorname{span}\lc \mathbf{v}_1,\mathbf{w}_1,\mathbf{v}_2,\mathbf{w}_2,\ldots,\mathbf{v}_m,\mathbf{w}_m \rc $ and $\operatorname{dim}(\operatorname{ker}P(D))=2m$.

\end{proof}

\begin{corollary}
\label{DsquaredplusI}
Consider the differential operator $Q=\lp\frac{L}{2\pi}\rp^2\ddxk{2}+m^2$ where $m\in \{ 1,2,\ldots,n \}$. The representation $\lp\frac{L}{2\pi}\rp^2D^2+m^2I$ of $Q$ in $\mathbb{R}^{2n+1}$ has $$\operatorname{rank}\lp \lp\frac{L}{2\pi}\rp^2D^2+m^2I\rp=2n-1,$$
where $D$ is the $(2n+1)$-dimensional matrix that represents $\ddx{}$ in $\mathbb{R}^{2n+1}$ and $I=I_{2n+1}$.
\begin{proof}
Using theorem \ref{DtotheN} with $P(x)=\lp\frac{L}{2\pi}\rp^2x^2+m^2$, we compute $\phi(k)=m^2-k^2$ and $\zeta(k)=0$.  System \eqref{phizeta} has exactly one solution, $k=m$, in the set $\{ 1,2,\ldots,n\}$.  Thus, $\operatorname{rank}\lp \lp\frac{L}{2\pi}\rp^2 D^2+m^2 I \rp=2n-1$.
\end{proof}
\end{corollary}

\begin{corollary}
\label{rankDtothep}
Let $D$ be a $(2n+1)$-dimensional square matrix representing $\ddx{}$ in $\mathbb{R}^{2n+1}$. Then $\operatorname{rank}(D^p)=2n$ for all $p \in \mathbb{N}$.
\begin{proof}
The result follows from the application of theorem \ref{DtotheN} to the polynomial $P(x)=x^p$.  Note that $\alpha_0=0$, one among $\lc \phi(k),\zeta(k) \rc$ equals $\pm \lp \frac{2\pi}{L} k \rp^p$, and the other equals zero.  Because system \eqref{phizeta} has no solutions in the set $\lc 1,2,\ldots,n \rc$, $\operatorname{rank}P(D)=2n$.
\end{proof}
\end{corollary}

\section{Numerical tests}
\subsection{Boundary value problems}
\begin{exmp}
\label{michaelsExample}
The aim of this example is to test accuracy and convergence rate of the trigonometric interpolation method.
Consider the boundary value problem 
\begin{equation}
\label{Example1}
\lc \begin{array}{l} u''+u=0,\\
u(0)=u(\pi/2)=1, \end{array}\right.
\end{equation} for a function $u:\lb 0,\frac{\pi}{2} \rb \rightarrow \mathbb{R}$. Let $\sigma = \lc x_i \rc_{i=0}^{2n}$ be a partition of $\lb 0,2\pi \rp$ such that $x_k=\frac{\pi}{2}$ for some $k\in\lc 0,1,\ldots,2n\rc$ and let $D, I \in \R^{(2n+1) \times (2n+1)}$ be the finite-dimensional representations of operators $\ddx{}$ and $1,$ constructed using this partition and the period $L=2\pi$. We approximate the differential operator $Q=\ddxk{2}+1$ with domain $C^2 [0,2\pi]$ by $D^2+I$ and reduce problem \eqref{Example1} to the problem of finding a vector $\hat{\mathbf{u}} \in\mathbb{R}^{2n+1}$ that satisfies
\begin{equation}
\label{DSquaredI}
\lc \begin{array}{l}
(D^2+I)\hat{\mathbf{u}}=\mathbf{0},\\
\hat{\mathbf{u}}_0=\hat{\mathbf{u}}_{k}=1.  
\end{array}\right.
\end{equation}

By Corollary \ref{DsquaredplusI}, the matrix $D^2+I$ has rank $2n-1$, so the space of solutions of the first equation in \eqref{DSquaredI} is two-dimensional, which implies the uniqueness of the solution of system~\eqref{DSquaredI}. To solve system \eqref{DSquaredI} in MATLAB, we define a new $(2n+3)\times(2n+1)$ matrix C whose upper block is $D^2+I$ and lower block consists of two rows $\mathbf{e}_1, \mathbf{e}_k \in \mathbb{R}^{2n+1}$, where $\{\mathbf{e}_i\}_{i=1}^{2n+1}$ is the standard basis of $\R^{2n+1}$. Let $\hat{\mathbf{f}}\in\mathbb{R}^{2n+3}$ be defined as $\hat{\mathbf{f}}=\lp 0, 0, \ldots, 0, 1, 1 \rp^T$.
We are now looking for the solution $\hat{\mathbf{u}}\in\mathbb{R}^{2n+1}$ to a system of linear equations $C\hat{\mathbf{u}}=\hat{\mathbf{f}}.$ We used MATLAB to solve this system and compared the solution $\hat{\mathbf{u}}$  with the vector of values of the exact solution $u(x)=\cos{x}+\sin{x}$ of \eqref{Example1} at the $N=2n+1$ points of the partition $\sigma$. The values of the error
$$E_{\max}(N)=\displaystyle\max\limits_{j\in \lc 0,1,\ldots,2n\rc} \left|\hat{\mathbf{u}}_j- u\lp x_j\rp \right|$$
are given in table \ref{ErrorTable}.

\begin{table}[h]
\caption{The errors $E_{\max}(N)$ obtained using trigonometric interpolation, Lagrangian interpolation, and the shooting method for boundary value problem \eqref{Example1}, where $N$ is the number of equidistant partition points of $\lb 0, 2\pi \rp$.}
\begin{tabular}{c|c|c|c|c}
\label{ErrorTable}
Method   & N& 5 & 7 & 9  \\ \hline
Trigonometric Interp. & $E_{\max}(N)$ & 6.9389e-16 & 3.7748e-15 & 3.9968e-15  \\
Lagrangian Interp. & $E_{\max}(N)$ & 1.1466e-004 & 5.2463e-006 & 2.5575e-007\\
Shooting Method & $E_{\max}(N)$ & 7.3e-3 & 3.2e-3 & 1.8e-3

\end{tabular}
\end{table}
We compared the accuracy of the trigonometric interpolation method described in this paper with the accuracy of the Lagrange interpolation method~\cite{Bihun,Calogero} and the standard shooting method.  To apply the shooting technique, we recast boundary value problem \eqref{Example1} into two initial value problems:  
\begin{equation*}
\lc \begin{array}{l}
v''=-v, \\ v(0)=1, v'(0)=0
\end{array} \right.
\end{equation*}
and 
\begin{equation*}
\lc \begin{array}{l}
w''=-w, \\ w(0)=0, w'(0)=1,
\end{array} \right.
\end{equation*}
which we solve by the method of finite differences.  The solution to problem \eqref{Example1} can be found by the formula

\begin{equation*}
u=v+\frac{1-v(\frac{\pi}{2})}{w(\frac{\pi}{2})}w.
\end{equation*}
Trigonometric interpolation is a suitable choice for problem \eqref{Example1}, since its solution is $2\pi$-periodic. In fact, $u(x_i)=\hat{\mathbf{u}}_i$ for all $i \in \{0,1,\ldots,n\}$ if $\hat{\mathbf{u}}$ is the \textit{exact} solution of $C\hat{\mathbf{u}}=\hat{\mathbf{f}}$ and the errors in the first row of table~\ref{ErrorTable} are computational errors that 
occur due to the limited number (sixteen) of significant digits of each component of $\hat{\mathbf{u}}$ stored in the computer memory.
Table~\ref{ErrorTable} shows that in this case our method yields significantly higher accuracy when compared to the shooting method or a similar method that uses Lagrangian interpolation \cite{Bihun,Calogero}.

\end{exmp}

\begin{exmp}
\label{uDoublePrimeProblem}
Consider the boundary value problem 
\begin{equation}
\label{uuuproblem}
\lc \begin{array}{l} u''+u'+u=0,\\
u(0)=u(\frac{\pi}{2\sqrt{3}})=1. \end{array}\right.
\end{equation}
Let $D$ be the finite representation of $\ddx{}$ given by formula~\eqref{D}, where $L=\frac{\pi}{\sqrt{3}}$.  By Corollary~\ref{DsquaredplusI}, the rank of $D^2+D+I$ is full, so we make the substitution $$u(x)=x \lp x-\frac{\pi}{2\sqrt{3}} \rp v(x)+1$$ to eliminate the boundary conditions in \eqref{uuuproblem}. The substitution reduces problem~\eqref{uuuproblem} to the differential equation
\begin{equation}
\label{uuusubproblem}
[(2 -\alpha )(x+ 1)+x^2]v+ \lp -2\alpha +4x-\alpha x+x^2 \rp v'
+\lp -\alpha x +x^2 \rp v'' = -1,
\end{equation}
where $\alpha= \frac{\pi}{2\sqrt{3}}$, with the conditions that require $v(0)$ and $v(\alpha)$ to be finite.

To approximate problem~\eqref{uuusubproblem} with a system of linear equations using trigonometric interpolation, we must choose the period $L$ of the trigonometric functions used in the approximation scheme. Although a natural choice seems to be $L=\frac{\pi}{2\sqrt{3}}$, it would impose an implicit condition $v(0)=v\lp \frac{\pi}{2\sqrt{3}} \rp$. Having this in mind, we double the length of approximation interval and use the period $L=\frac{\pi}{\sqrt{3}}$ in formula~\eqref{D}. We partition the interval $\lb 0, \frac{\pi}{2\sqrt{3}} \rb,$ into $(2n+1)$~points $x_j=\frac{\pi}{2\sqrt{3}}\cdot \frac{j}{2n}$, where $j\in \{0,1,\ldots,2n\}$, and consider this as a partition of $[ 0, \frac{\pi}{\sqrt{3}})$. We compute the matrix ${D}$
using formula~\eqref{D}. To improve the conditioning of a matrix representation of $\lp \ddxk{2}+\ddx{}+1\rp$, we use the matrix $\hat{D}=\Psi^{-1} D \Psi$ instead of ${D}$, where $\Psi=\operatorname{Diag}({\psi}_0,{\psi}_1,\ldots,{\psi}_{2n})$ (read Theorem~\ref{dhatworks} for more insight on the use of $\hat{D}$).

 We then solve the system of linear equations
\begin{align*}
&[(2 -\alpha)(X+ I)+X^2+ \lp -2\alpha I+4X-\alpha X+X^2 \rp \hat{D}\\ 
&+\lp -\alpha X +X^2 \rp \hat{D}^2 ]\hat{\mathbf{v}} = -\operatorname{diag} I
\end{align*} 
for $\hat{\mathbf{v}} \in \mathbb{R}^{2n+1}$, and calculate the vector $\hat{\mathbf{u}}\in \R^{2n+1}$ using formula $\hat{\mathbf{u}}_i=x_i\lp x_i -\frac{\pi}{2\sqrt{3}} \rp{\psi}_i\hat{\mathbf{v}}_i+1$, where $i \in \{0,1,\ldots,2n\}$. The vector $\hat{\mathbf{u}}$ approximates the values of the exact solution
$$ u(x)=\exp\lp{-\frac{x}{2}}\rp \lb  \lp \sqrt{2}\exp{\frac{\pi}{4\sqrt{3}}}-1 \rp 
\sin  \frac{ \sqrt{3} x }{2}   
 +\cos \frac{\sqrt{3}x}{2} \rb.$$
 to problem~\eqref{uuuproblem} at the partition points.

In the following table  we summarize the approximation error $$E_{\max}(N)=\max_{i\in \{0,1,\ldots,2n\}}\left|\hat{\mathbf{u}}_i-u(x_i)\right|$$ for several values of $N=2n+1$. 

\begin{table}[h]
\caption{Values of the error $E_{\max}$ obtained when solving problem~\eqref{uuuproblem} using trigonometric interpolation.}
\begin{tabular}{c|c|c|c}
\label{ErrorTable2}
 Number of Partition Points N& 5 & 11 & 21  \\ \hline
 $E_{\max}(N)$ & 5.6458e-4 & 3.7868e-6 & 3.7706e-9  
\end{tabular}
\end{table}

\end{exmp}

\begin{exmp}
\label{noSubExample1}
Consider the boundary value problem
\begin{equation}
\label{nosubproblem1}
\lc \begin{array}{l} \lp \frac{1}{2\pi} \rp^4 u^{(4)} + 13 \lp \frac{1}{2\pi} \rp^2 u^{(2)}  + 36 u=0,\\
u(0)=u(1)=1, u'(0)=u'(1)=0. \end{array}\right.
\end{equation}

We partition the interval $[0,1)$ by points $\{x_j\}_{j=0}^{2n}$ and approximate the differential equation in \eqref{nosubproblem1} by the system of linear equations $A\hat{\mathbf{u}}=\mathbf{0}$ for $\hat{\mathbf{u}} \in \mathbb{R}^{2n+1}$, where 
$$A=\lp \frac{1}{2\pi} \rp^4 {D}^4+ 13 \lp \frac{1}{2\pi} \rp^2 {D}^2 + 36I $$ 
and ${D}$ is computed by formula~\eqref{D} using the period $L=1$.

By theorem \ref{DtotheN}, $\operatorname{rank}A=\operatorname{order}A-4$. Therefore, we can complement the linear system $A\hat{\mathbf{u}}=\mathbf{0}$ with four additional equations that represent the boundary conditions in \eqref{nosubproblem1}. More precisely, we solve the linear system
$$
\lc \begin{array}{l} A\hat{\mathbf{u}}=\mathbf{0}, \\ 
\hat{\mathbf{u}}_0=\hat{\mathbf{u}}_{2n}=1, \\ 
\lb{D}\hat{\mathbf{u}} \rb_0=\lb{D}\hat{\mathbf{u}} \rb_{2n}=0 
 \end{array}\right.
$$
or, equivalently, 
$$\lp
\begin{array}{c}
A \\ 
B \\ 
\end{array}\rp 
\hat{\mathbf{u}}=
\lp
\begin{array}{c}
\mathbf{0} \\
\mathbf{b}
\end{array}\rp,
$$

where the block 

$$
B=\lp
\begin{array}{ccccc}
1 & 0 & \ldots & 0 & 0 \\
0 & 0 & \ldots & 0 & 1 \\
{D}_{0,0} &{D}_{0,1} & \ldots &{D}_{0,2n-1} &{D}_{0,2n} \\
{D}_{2n,0} &{D}_{2n,1} & \ldots &{D}_{2n,2n-1} &{D}_{2n,2n} \\
\end{array}\rp
$$

and $\mathbf{b}=\lp 1,1,0,0\rp^T$. 

We compute the approximate solution of \eqref{nosubproblem1} by the formula 
\begin{equation}
\label{approxfunction}
u_{\operatorname{approx}}(x)= \sum_{i=0}^{2n} \hat{\mathbf{u}}_i{t}_i(x).
\end{equation}

In fig.~\ref{ApproximateSolutionNoSub1Figure}, we graph the approximate solution $u_{\operatorname{approx}}$ obtained using the partition $x_i=\frac{i}{2n+1}$, where $i \in \{ 0, 1, \ldots, 2n\}$ and $n=5$.

\begin{figure}[h]
\caption{The approximate solution to \eqref{nosubproblem1}.}
\label{ApproximateSolutionNoSub1Figure}
\includegraphics[scale=.75]{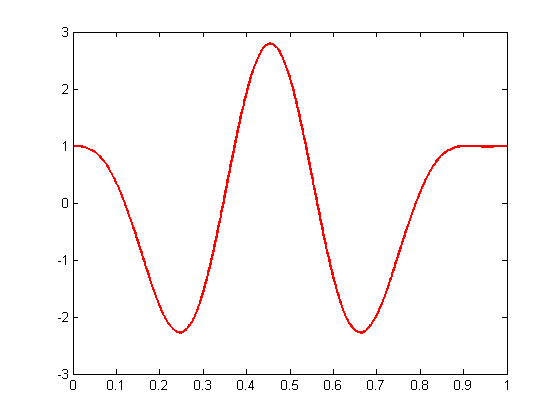}
\end{figure}   
\end{exmp}

\subsection{Eigenvalue problems}
\label{eigensubsection}
Eigenvalue problems play an important role in Physics and Applied Mathematics. Solutions to many partial differential equations that model heat conduction, wave propagation, and other physical phenomena, can be reduced to the solution of a Sturm-Liouville eigenvalue problem.  

In this section, we present solutions to several Sturm-Liouville eigenvalue problems obtained using trigonometric interpolation. Below we state some background information useful for understanding applications.

\begin{theorem}
\label{calogero}
Let $\mathcal{A}=\alpha_0+\alpha_1\ddx{}+\ldots+\alpha_s\ddxk{s}$ be a differential operator with domain $D(\mathcal{A}) \supset 
C^s_{[0,L]}$, where $\alpha_0,\alpha_1,\ldots,\alpha_s \in \R$ and $\alpha_s \neq 0$. Let $\lp \alpha,f\rp \in \mathbb{R}\times D(\mathcal{A})$ be an eigenpair of the operator $\mathcal{A}$. If $f \in \mathcal{T}_n^L$, then the finite-dimensional representation $A=a_0 + a_1D+\ldots+a_sD^s$ of operator $\mathcal{A}$, constructed using a partition $\{ x_i\}_{i=0}^{2n}$ of $[0,L)$, has an eigenpair $(\alpha,\mathbf{w})\in \mathbb{R}\times\mathbb{R}^{2n+1}$ such that $\sum_{i=0}^{2n}\mathbf{w}_it_i(x)=f(x)$ for all $x \in [0,L)$.
\end{theorem}

This theorem was proven by Calogero~\cite{Calogero3}. Because we use different notation, we supply a short proof in appendix~\ref{appb} for the reader's convenience.  

Let $\Psi=\operatorname{Diag}(\psi_0,\psi_1,\ldots,\psi_{2n})$, where $\psi_j$'s are defined by formula~\eqref{piconstant}.  In the following examples, we use the matrix $\hat{D}=\Psi^{-1}D\Psi$ instead of $D$ to solve eigenvalue problems because it yields representations of differential operators that are usually better conditioned. 

\begin{theorem}
\label{dhatworks}
If $\lp \alpha,f\rp \in \mathbb{R}\times \mathcal{T}_n^L$ is an eigenpair of a differential operator $\mathcal{A}$ defined as in theorem \ref{calogero}, then the matrix $\hat{A}=a_0+a_1\hat{D}+\ldots+a_s\hat{D}^s$ has an eigenpair $(\alpha,\mathbf{w})\in \mathbb{R}\times\mathbb{R}^{2n+1}$, where $\summ{i}{0}{2n} \psi_i \mathbf{w}_i t_i(x) = f(x)$ for all $x\in [0,L).$
\end{theorem}

This theorem is a direct consequence of Theorem~\ref{calogero} and the formula $D=\Psi \hat{D} \Psi^{-1}$.  

%Eigenvalueproblem1 in matlab%
\begin{exmp}
\label{EigenEXAMPLE1} This example aims to test the accuracy and convergence rate of the trigonometric interpolation method applied to eigenvalue problems.  
Consider the eigenvalue problem 
\begin{equation}
\label{EigExample1}
\lc \begin{array}{l} -y''=\lambda y,\\
y(0)=y(1)=0. \end{array}\right.
\end{equation} Let $\{x_i\}^{2n}_{i=0}$ be a partition of $[0,1)$. Let the matrix $\hat{D}=\Psi^{-1}D\Psi$, where $D$ is defined by eq.~\eqref{D} with $L=1$. Let $X=\operatorname{Diag}\lp x_0,\ldots,x_{2n} \rp$ and $I=I_{2n+1}$. To get rid of boundary conditions in \eqref{EigExample1}, we make a substitution $y(x)=x(x-1)v(x)$.
Via this substitution, the eigenvalue problem \eqref{EigExample1} is recast as 
$$-\left[(x^2-x)v''+(4x-2)v'+2v \right] = \lambda x(x-1)v,$$
with the condition that requires $v(0)$ and $v(1)$ to be finite. The last eigenvalue problem is then approximated by the matrix eigenvalue problem 
\begin{align}
\label{Aex_1}
-\lp (X^2-X)\hat{D}^2+(4X-2I)\hat{D}+2I \rp \hat{\mathbf{v}}=\lambda \lp X(X-I)\rp \hat{\mathbf{v}}, \qquad \hat{\mathbf{v}} \in \R^{2n+1}.  
\end{align}
We use MATLAB to find the eigenvalues of \eqref{Aex_1} using the partition $x_i=\frac{i}{n}$ for $i \in \{0,1,\ldots,n-1\}$. The exact eigenvalues for this system are known to be $\pi^2j^2$ for $j\in\mathbb{N}$, and we compare the approximate eigenvalues with exact eigenvalues in table \ref{ErrorTable2}, which also illustrates the convergence of this scheme as the number of equidistant partition points increases.

\begin{table}[h]
\caption{Approximate eigenvalues of problem \eqref{EigExample1}, where $\lambda_i$ denotes the $i^{th}$ eigenvalue.}
\begin{tabular}{c|c|c|c|c}
\label{ErrorTable2}
Number of partition points   & $\lambda_1$ & $\lambda_5$ & $\lambda_{11}$ & $\lambda_{15}$  \\ \hline
19 & 9.906 & 247.699 & 1200.521 & 2240.381  \\
31 & *.883 & **7.088 & *196.079 & **24.557  \\
51 & *.*75 & **6.867 & ***4.858 & ***1.891  \\
101 & *.**1 & ***.772 & ****.38 & ***0.956  \\ \hline
Exact eigenvalue & 9.8696 & 246.740 & 1194.222 & 2220.660  
\end{tabular}
\end{table}
\end{exmp}
Using the trigonometric interpolation approach, we solved several eigenvalue problems and compared the approximate eigenvalues we obtained with the lower and upper bounds found in \cite{Fichera}.

\begin{exmp}
\label{EigenEXAMPLE2}
Consider a genetics eigenvalue problem 
\begin{equation}
\label{eqEigEx2}
\lc \begin{array}{l} 
\frac {1}{4} \lp 1-x^2 \rp \lp -u''
+2\alpha \lb \lp x_0 - x \rp u \rb' \rp = \lambda u,\\
u(-1)=u(1)=0, \end{array}\right.
\end{equation}
where $ x_0, \alpha \in \mathbb{R}$.

Using the substitution $u(x)= (1-x^2)v(x),$ we approximate \eqref{eqEigEx2} by the matrix eigenvalue problem $A\hat{\mathbf{v}}=\lambda \hat{\mathbf{v}}$, where $\hat{\mathbf{v}} \in \R^{2n+1}$ and

\begin{align}
\label{Aex_2}
A&=-\frac{1}{4}\lp I-X^2 \rp \hat{D}^2 + \lb X+ \frac{1}{2}\alpha \lp I-X^2\rp \lp x_0 I - X\rp \rb \hat{D}  \\
&+ \frac{1}{2} + \alpha X \lp X-x_0 I \rp - \frac{1}{2} \alpha \lp I-X^2 \rp. \nonumber
\end{align}
Results in tables \ref{ErrorTableEx3151} and \ref{ErrorTableEx3152} are obtained using the partition  $x_i=\frac{2i-n-1}{n}$ for all $i \in \{ 1,2,\ldots,n \}$, where $n=25$ and the period $L=4$.

\begin{table}
\caption{$\lambda_i$ denotes the $i^{th}$ eigenvalue of problem \eqref{eqEigEx2}, where $\alpha= 1$ and $x_0=0.6$. Upper and lower bounds for $\lambda_i$ are taken from \cite{Fichera}. }
\begin{tabular}{c|c|c|c}
\label{ErrorTableEx3151}
   & Approximate Eigenvalues & Lower Bound & Upper Bound \\ \hline
$\lambda_1$ & 0.3938 & 0.3931 & 0.3939  \\
$\lambda_4$ & 4.9514 & 4.85 & 4.96  \\
$\lambda_8$ & 17.9515 & 16.7 & 17.96  \\
$\lambda_{10}$ & 27.4552 & 24.6 & 27.46     
\end{tabular}
\end{table}

\begin{table}
\caption{$\lambda_i$ denotes the $i^{th}$ eigenvalue of problem \eqref{eqEigEx2}, where $\alpha= 2$, and $x_0=0.4$. Upper and lower bounds for $\lambda_i$ are taken from \cite{Fichera}.  }
\begin{tabular}{c|c|c|c}
\label{ErrorTableEx3152}
   & Approximate Eigenvalues & Lower Bound & Upper Bound \\ \hline
$\lambda_1$ & 0.2991 & 0.2982 & 0.2992  \\
$\lambda_4$ & 4.9625 & 4.73 & 4.97  \\
$\lambda_8$ & 17.9570 & 15.3 & 17.96  \\
$\lambda_{10}$ & 27.4585 & 21.7 & 27.46     
\end{tabular}
\end{table}
\end{exmp}

\begin{exmp}
\label{EigenEXAMPLE3}
Consider the eigenvalue problem 
\begin{equation}
\label{eigEquEx3}
\lc \begin{array}{l} u''+\lambda \lp 1+\sin x \rp u = 0,\\
u(0)=u(\pi)=0. \end{array}\right.
\end{equation}
We make the substitution $u(x)=y(1-y)v(y)$, where $y=x/\pi$. We approximate problem~\eqref{eigEquEx3} by the matrix eigenvalue problem $A\hat{\mathbf{v}}=\lambda \hat{\mathbf{v}}$, where $\hat{\mathbf{v}} \in \R^{2n+1}$ and
\begin{align*}
\label{Aex_3}
A = \frac{1}{\pi^2}\lp I+\sin(\pi X)\rp^{-1} \lp -\hat{D}^2+2 \lb X \lp I-X \rp\rb^{-1}\lb \lp 2X-I\rp \hat{D} +1 \rb \rp.
\end{align*}
The matrices $X$ and $\hat{D}$ are constructed using a partition $\{x_j\}_{j=0}^{2n}$ of $[0,2)$. We use the interval $[0,2)$ instead of $[0,1)$ for approximation of the function $v$ to avoid imposing the implicit condition $v(0)=v(1)$.

Results in table~\ref{ErrorTableEx325} are obtained using the partition $x_j=\frac{j}{n+1}$ for $j \in \{ 1,2,\ldots,n \}$, $n=21$, and the period $L=2$.

\begin{table}
\caption{$\lambda_i$ denotes the $i^{th}$ eigenvalue of problem \ref{eigEquEx3}. Upper and lower bounds for $\lambda_i$ are taken from \cite{Fichera}.}
\begin{tabular}{c|c|c|c}
\label{ErrorTableEx325}
   & Approximate Eigenvalues & Lower Bound & Upper Bound \\ \hline
$\lambda_1$ & 0.5403 & 0.540282 & 0.540319  \\
$\lambda_3$ & 5.4486 & 5.411 & 5.449  \\
$\lambda_5$ & 15.3126 & 14.54 & 15.313  \\
$\lambda_{6}$ & 22.0972 & 21.39 & 22.1     
\end{tabular}
\end{table}
\end{exmp}

\section{Acknowledgments}

This project is supported by the Centennial Scholars Research Grant of Concordia College,~MN, PI~Oksana Bihun, and the NSF grant DUE~0969568, PI~Heidi Manning.

\appendix

\section{Proof of theorem~\ref{calogero}}
\label{appb}
\begin{proof}
Let $(\alpha,f)\in \mathbb{R}\times \mathcal{T}_n^L$ be an eigenpair of $\mathcal{A}.$ Choose a partition $\{x_i\}_{i=0}^{2n}$ of $[0,L)$ and let 
$A=a_0+a_1D+\ldots+a_sD^s, $ 
where $D$ is given by formula~\eqref{D}.
Let $\hat{\mathbf{f}}=\lp f(x_0),f(x_1),\ldots,f(x_{2n})\rp^T$.

Because $f \in \mathcal{T}_n^L$, $\mathcal{A}f \in \mathcal{T}_n^L[x]$ and the equation $\mathcal{A}f(x)=\alpha f(x)$ is equivalent to the equation $$\sum_{j=0}^{2n} \lb A\hat{\mathbf{f}}\rb_j t_j(x)=\sum_{j=0}^{2n} \alpha \hat{\mathbf{f}}_jt_j(x), \qquad x \in [0,L),$$ 
by corollary~\ref{TofFisF}. The last equation is equivalent to $A\hat{\mathbf{f}}=\alpha \hat{\mathbf{f}}.$  Hence, $(\alpha,\hat{\mathbf{f}})\in \mathbb{R} \times \mathbb{R}^{2n+1}$ is an eigenpair of $A$. Moreover, $f(x)=\sum_{i=0}^{2n}\hat{\mathbf{f}}_i t_i(x)$ for all $ x \in [0,L)$ as required.
\end{proof}

\end{document}